\documentclass[12pt]{amsart}
\usepackage{amssymb}
\usepackage{amsmath}
\usepackage{amsfonts}
\usepackage{latexsym}
\usepackage[T1]{fontenc}
\usepackage[latin1]{inputenc}
\usepackage{hyperref}
\usepackage[all,ps]{xy}
\RequirePackage{amsthm}
\RequirePackage{amssymb}
\usepackage[all,ps]{xy}
\usepackage{latexsym}
\usepackage{amscd}
\newcommand{\Pic}{{\rm Pic}}
\renewcommand{\lg}{{\rm lg}}
\newcommand{\gon}{{\rm gon}}
\newcommand{\Cliff}{{\rm Cliff}}
\newcommand{\supp}{{\rm supp}}
\newcommand{\Gr}{{\rm Gr}}
\newcommand{\rk}{{\rm rk}}

\renewcommand{\det}{{\rm det}}
\renewcommand{\Im}{{\rm Im}}
\newcommand{\Ann}{{\rm Ann}}
\renewcommand{\dim}{{\rm dim}}
\renewcommand{\deg}{{\rm deg}}
\newcommand{\sO}{{\mathcal O}}
\newcommand{\sE}{{\mathcal E}}

\newcommand{\ev}{{\rm ev}}

\begin{document}

\newtheorem{thm}{Theorem}[section]
\newtheorem{defn}[thm]{Definition}
\newtheorem{subthm}[thm]{Sub-Theorem}
\newtheorem{warning}[thm]{Warning}
\newtheorem{summary}[thm]{Summary}
\numberwithin{equation}{subsection}
\newtheorem{claim}[thm]{Claim}
\newtheorem{cor}[thm]{Corollary}
\newtheorem{rmk}[thm]{Remark}
\newtheorem{prop}[thm]{Proposition}
\newtheorem{lem}[thm]{Lemma}
\newtheorem{sublem}[equation]{Lemma}
\newtheorem{ex}[thm]{Example}
\newtheorem{con}[thm]{Conjecture}

\title{The Green Conjecture for Exceptional Curves on a $K3$ Surface}
\author
{Marian Aprodu (IMAR) and Gianluca Pacienza (IRMA)}
\date{\today}
\thanks{MA was supported in part by a Humboldt Return
Fellowship and by the ANCS grant 2-CEx 06-11-20/2006.}
\maketitle

\begin{abstract}
We use the Brill-Noether theory to prove the Green conjecture for
exceptional curves on $K3$ surfaces. Such curves count among the few
ones having Clifford dimension $\ge 3$. We obtain our result by
adopting an infinitesimal approach due to Pareschi, and using the
degenerate version of the Hirschowitz-Ramanan-Voisin theorem
obtained in \cite{A05}.
\end{abstract}

\section{Introduction.}

Two conjectures made in the eighties by Green, and Green-Lazarsfeld,
pointed out to some deep links between the intrinsic properties of
algebraic curves, and their Koszul cohomology groups with values in
suitably chosen line bundles. Recall that the Koszul cohomology
$K_{p,q}(X,L)$ of a complex projective variety $X$ with values in a
globally generated line bundle $L$ on $X$ is defined as the
cohomology at the middle of the complex:

$$
\bigwedge^{p+1}H^0(L)\otimes H^0(L^{\otimes(q-1)})\to
\bigwedge^pH^0(L)\otimes H^0(L^{\otimes q})\to
\bigwedge^{p-1}H^0(L)\otimes H^0(L^{\otimes(q+1)}).
$$

A basic result due to Green and Lazarsfeld, see
\cite[Appendix]{G84}, shows that $K_{r_1+r_2-1,1}(X,L)$ is not zero
if $X$ is smooth and $L$ is decomposed as $L=L_1+L_2$ with
$r_i:=h^0(X,L_i)-1\ge 1$. In the particular case of curves, one
obtains that $K_{g-c-2,1}(X,K_X)$ and $K_{h^0(L)-d-1,1}(X,L)$ are
not zero, where $g$ is the genus of $X$, $c$ is the Clifford index
(see \cite{Ma82} for the definition), $d$ is the gonality, and $L$
is an arbitrary line bundle of sufficiently large degree. The
above-quoted conjectures predict that these non-vanishing results
are sharp.

\begin{con}[Green, \cite{G84}]
$K_{g-c-1,1}(X,K_X)=0$
\end{con}

\begin{con}[Green-Lazarsfeld, \cite{GL85}]
$K_{h^0(L)-d,1}(X,L)=0$
\end{con}

In recent years, strong evidence has been found for these
conjectures. Putting together Voisin's \cite{V02} and Teixidor's
\cite{T02} results, Green conjecture is valid for a generic
$d$-gonal curve, for $d<[g/2]+2$. The case of generic curves of odd
genus, not covered by the above-mentioned results was settled by
Voisin \cite{V05}. This case was particularly challenging due to a
previous work of Hirschowitz and Ramanan \cite{HR98}, which together
with Voisin's result \cite{V05} implies the validity of the Green
conjecture for {\it any} curve of odd genus and maximal gonality.
The Green-Lazarsfeld conjecture is also valid for generic $d$-gonal
curves, see \cite{AV03}, \cite{A04}, and \cite{A05}.

\medskip

The Brill-Noether theory could play a major r\^ole in the attempt to
solve these conjectures. Specifically, suppose $d\ge 3$ is an
integer, and $X$ is a smooth $d$-gonal curve with $d<[g/2]+2$, such
that the varieties $W^1_{d+n}(X)$, parametrising degree-$(d+n)$
pencils on $X$,  verify

\medskip

\noindent {\bf{(Linear growth conditions):}}
$\dim (W^1_{d+n}(X)) \leq n,$ for all $n$ such that
$ 0\le n\le g-2d+2.$

\medskip

Then $X$ verifies both Green, and Green-Lazarsfeld conjectures, see
\cite[Theorem 2]{A05}.

\medskip

It is therefore important to try and control the dimensions of the
varieties of pencils, and find intervals on which their growth is
linear in the degree. A prediction was made by G. Martens,
\cite[Statement (T)]{Ma84}:

\medskip

\noindent{\bf Statement (T).} {\it If $j\ge 1$, $g\ge 2j+4$, and
$m\in[j+3,g-1-j]$ are integer numbers, and $X$ is a curve of genus
$g$, such that $\dim(W^1_{m}(X))=m-2-j$, then
$\dim(W^1_{s}(X))=s-2-j$ for all integers $s\in [j+3,g-j]$.}

\medskip

Notice that one of the conditions above is $\dim(W^1_{j+3}(X))=1$,
and the inequality $d\ge j+3$ yields to $j+3\le (g+1)/2$. In
particular, for Brill-Noether generic curves, the Statement (T) is
empty. For explicit special curves of Clifford dimension one,
Statement (T) with $m=\gon(X)+1$, together with \cite[Theorem
2]{A05} can be seen as a potential tool for verifying the Green, and
Green-Lazarsfeld conjectures for new classes of curves.

\medskip

Partial versions of Statement (T) have been proved (see for instance
\cite{Ho82}, \cite{Ma84}, \cite{CKM92}). Unfortunately, the picture
happens to be more subtle than expected (and hoped). We exhibit
counter-examples to Statement (T) as special curves in the linear
system $|2H|$ on a $K3$ surface whose Picard group is generated by a
hyperplane section $H$ and a line $\ell$.

\begin{prop}
\label{counter-example} Let $S$ be a $K3$ surface with
$\Pic(S)=\mathbb{Z}.H\oplus \mathbb{Z} \ell$, with $H$ very ample,
$H^2=2r-2$, and $H.\ell=1$. There exists a smooth curve $C\in
|2H|$ (whose genus equals $4r-3$, and whose gonality equals $2r-2$,
see \S 3.1), and
\begin{enumerate}
\item $\dim(W^1_{2r-2}(C))=1$.
\item $\dim(W^1_{2r-1}(C))=2$.
\item $\dim(W^1_{2r}(C))=4$.
\item $\dim(W^1_{2r+1}(C))$ equals $5$ or $6$.
\end{enumerate}
\end{prop}

Nevertheless, we prove that the generic curves in $|2H|$ verify a
weak linear growth condition.

\begin{thm}
\label{2h} Let $S$ be a $K3$ surface with
$\Pic(S)=\mathbb{Z}.H\oplus \mathbb{Z} \ell$, with $H$ very ample,
$H^2=2r-2\ge 4$, and $H.\ell=1$. Then the generic curve in the
linear system $|2H|$ verifies
$$
\dim(W^1_{2r-2+n}(C))=n,\ \ \ {\textrm {for}}\ n\in\{0,1,2\}.
$$
\end{thm}

A consequence of this result is that the Green-Lazarsfeld conjecture
holds for the generic curves in the linear system $|2H|$, see
Corollary \ref{GLfor2H}.

The idea of the proof is to look at the family of pairs $(C,A)$,
with $C\in |2H|$ smooth and $A\in W^1_{2r-2+n}(C)$, and to give a
bound on the dimension of the irreducible components dominating
$|2H|$. Thanks to the work of Lazarsfeld and Mukai, to the data
$(C,A)$ (for simplicity we assume here that $A$ is a complete and
base-point-free pencil) one can attach a rank $2$ vector bundle
$E(C,A)$ on the surface $S$. If this bundle is {\it simple}, then
the original argument of Lazarsfeld's [L], or the variant provided
by Pareschi \cite{P95}, allows one to determine these dimensions. In
the {\it non-simple} case a
useful lemma (see \cite{GL87}, \cite{DM89} and \cite{CP95}), brings
to a very concrete description of the parameter space for such
bundles. This description, together with the infinitesimal
approach of Pareschi \cite{P95}, allows us to conclude.

Moreover, and more importantly, using the  degenerate version of the
Hirschowitz-Ramanan-Voisin theorem  \cite[Proposition 8]{A05}, we
derive from Theorem \ref{2h} the following.

\begin{thm}
\label{2h+l} Let $S$ be a $K3$ surface with
$\Pic(S)=\mathbb{Z}.H\oplus \mathbb{Z} \ell$, with $H$ very ample,
$H^2=2r-2\ge 4$, and $H.\ell=1$. Then any smooth curve in the linear
system $|2H+\ell|$ verifies the Green conjecture.
\end{thm}

Smooth curves in the linear system $|2H+\ell|$ are particularly
interesting in several regards. First of all, they count among the
few examples of curves whose Clifford index is not computed by
pencils, i.e. $\textrm{Cliff}(C)=\gon(C)-3$, as it was shown in
\cite{ELMS89} (other obvious examples are given by plane curves, for
which the Green conjecture is already checked, see \cite{Lo89}).
Such curves are the most special in the moduli space of curves
from the Clifford index view-point, reason
for which some authors call them {\it exceptional} curves. Hence,
the case of smooth curves in $|2H+\ell|$ may be considered as
opposite to that of generic curve of fixed gonality. Secondly, they
carry a one-parameter family of pencils of minimal degree (see
\cite{ELMS89}), so \cite[Theorem 2]{A05} cannot be applied directly.

The outline of the article is the following. First, we recall some
vector bundle techniques which will be used in the proof of Theorems
\ref{2h} and \ref{2h+l}. The proof of Theorem \ref{2h} splits in
several cases. In \S $3.1$ we treat the case $n=0,1$. In \S 3.2, we
analyse the case $n=2$, and $E(C,A)$ simple. The remaining case,
$n=2$, and $E(C,A)$ not simple are ruled out in \S 3.3. We put
together all these intermediate steps in \S 3.4 and prove Theorem
\ref{2h} and Proposition \ref{counter-example}. Finally, we show how
Theorem \ref{2h} implies  Theorem \ref{2h+l}.

\medskip

\noindent{\em Acknowledgments.} MA would like to thank I.R.M.A.
Strasbourg, and especially Olivier Debarre for the warm hospitality
during the first stage of this work. The authors are grateful to O.
Debarre for numerous conversations on this topic.

\section{Vector bundle techniques.}

In the present Section we recall some basic vector bundle
techniques, see \cite{La86}, \cite{GL87}, \cite{La89}, \cite{OSS80}.

\subsection{The Lazarsfeld-Mukai vector bundle.}

Given a smooth curve $C$ belonging to a linear system $|L|$ on a
$K3$ surface $S$ and a  base-point-free line bundle $A\in \Pic(C)$,
we recall  how to associate to this data, following \cite{La86} and
\cite{Mu89}, a vector bundle $E:=E(C,A)$ of rank $h^0(C,A)$ on $S$
and we  record a number of properties of $E$ which will be freely
used in the rest of the paper. First one considers the rank
$h^0(C,A)$ vector bundle $F(C,A)$ defined as the kernel of the
evaluation of sections of $A$ (considered as torsion sheaf on the
whole surface $S$)

\begin{equation}\label{F}
 0\to F(C,A)\to H^0(C,A)\otimes \sO_S \buildrel{ev}\over{\to} A\to 0.
\end{equation}

Then dualizing the above exact sequence and setting $E:=
E(C,A):=F(C,A)^*$ we get :

\begin{equation}
\label{E}
 0\to H^0(C,A)^*\otimes \sO_S \to E\to K_C(-A)\to 0.
\end{equation}

The invariants of $E$ are :
\begin{enumerate}
 \item $\det (E)= \sO_S(C)$;
 \item $c_2(E)=\deg(A)$;
 \item $h^0 (S,E)= h^0(C,A)+h^1(C,A)= 2h^0(C,A)-\deg(A)-1 +g(C).$
\end{enumerate}

Moreover, $E$ is globally generated off a finite set, and
$$
 h^1(S,E)=h^2(S,E)=0.
$$

There is a natural rational map
\begin{equation}\label{d_E}
 d_E : \Gr(\rk(E),H^0(S,E)) \dashrightarrow |L|
\end{equation}
from the Grassmanniann of $\rk(E)$-dimensional subspaces of
$H^0(S,E)$ to the linear system $|L|$. This map sends a generic
subspace $\Lambda\in \Gr(\rk(E),H^0(S,E))$ to the degeneracy locus
of the evaluation map :
$$
 \ev_{\Lambda} : \Lambda \otimes \sO_S \to E
$$
(notice that, generically, this degeneracy locus cannot be the whole
surface, since $E$ is generated off a finite set).

Moreover, for generic $\Lambda$, the image $d_E(\Lambda)$ is a
smooth curve $C_{\Lambda}$ (the smoothness of the degeneracy locus
is an open condition, and it is realized on a non-empty set of
$\Gr(\rk(E), H^0(S,E))$ by the construction of $E$), and the
cokernel of $\ev_{\Lambda}$ is a line bundle
$K_{C_\Lambda}(-A_{\Lambda})$ of $C_{\Lambda}$, where
$\deg(A_{\Lambda})=c_2(E)$. An important feature on the map $d_E$ is
that its differential at a point $\Lambda$ coincides with the
multiplication map
$$
 \mu_{0,A_{\Lambda}}: H^0(C_{\Lambda}, A_{\Lambda}) \otimes H^0(C_{\Lambda},
 K_{C_{\Lambda}} - A_{\Lambda})
 \to  H^0(C_{\Lambda}, K_{C_{\Lambda}}).
$$

On the other hand, if $M_A$ is the vector bundle (of rank
$h^0(C,A)-1$) on $C$ defined by the kernel of the evaluation map
{\it on the curve} :
\begin{equation}\label{MA}
 0\to M_A \to H^0(C,A)\otimes \sO_C \buildrel{ev}\over{\to} A\to 0,
\end{equation}
by tensoring (\ref{MA}) with $K_C(-A)$, one gets
\begin{equation}\label{kermu0&MA}
 \ker(\mu_{0,A})= H^0(C, M_A\otimes K_C(-A)).
\end{equation}

Notice that, by construction, there is a natural surjective map from
$F(C,A)_{|C}$ to $M_A$, and, by determinant reason, we have
\begin{equation}\label{FA&MA}
0\to A(-K_C)\to F(C,A)_{|C}\to M_A\to 0.
\end{equation}

\subsection{The Serre construction.}

The Serre construction consists in associating to a locally complete
intersection $0$-dimensional subscheme $\xi$ of $S$, and to a
non-zero element $t\in H^1(S,L\otimes I_\xi)^*$, where $L\in
\Pic(S)$, a rank two vector bundle $E:=E_{\xi, t}$ on $S$ and a
global section $s\in H^0(S,E)$, whose zero locus is $\xi$.

By the Grothendieck-Serre duality Theorem, we have
$$
H^1(S,L\otimes I_\xi)^*\cong\textrm{Ext}^1(L\otimes I_{\xi},\sO_S),
$$
hence to each  $t\in H^1(S,L\otimes I_\xi)^*$ we may associate a
sheaf $\mathcal E$ is given by an extension
$$
0\to \sO_S\to \mathcal E\to L\otimes I_\xi\to 0.
$$

Notice that the global section $s$ of $\mathcal E$ coming from the
inclusion $\sO_S\hookrightarrow \mathcal E$ vanishes on $\xi$.

The precise criterion for an extension as above to be locally free
is the following:

\begin{prop}[\cite{OSS80}, \cite{La97}]
\label{extension} Given $t\in H^1(S,L\otimes I_\xi)^*$, the
corresponding $\mathcal E$ fails to be locally free if and only if there
exists a proper (possibly empty!) subscheme $\xi'\subset \xi$ such
that
 $$
 t\in \Im(H^1(S,L\otimes I_{\xi'})^*\to H^1(S,L\otimes I_\xi)^*).
 $$
\end{prop}

We also recall a result which we state and only when $S$ is a $K3$
surface (note that the canonical bundle of $S$ is trivial).

\begin{thm}[\cite{GH78}]\label{GH}
There exists a rank two vector bundle $E$ on $S$, with $\det(E)=L$,
and a section $s\in H^0(S,E)$ such that $V(s)=\xi$ if, and only if,
every section of $L$ vanishing at all but one of the points in the
support of $\xi$ also vanishes at the remaining point.
\end{thm}

Theorem \ref{GH} immediately yields :

\begin{prop}\label{notreextension}
Let $S$ be a $K3$ surface, and $L$ a line bundle on $S$. Then, for
any $0$-dimensional subscheme $\xi$ of $S$, such that
$h^0(S,L\otimes I_{\xi'})=0$, for all $\xi'\subset \xi$ with
$\lg(\xi')=\lg(\xi)-1$, there exists a rank two vector bundle $E$ on
$S$ given by an extension
 $$
  0\to \sO_S\to E\to L\otimes I_\xi\to 0.
 $$
\end{prop}

\subsection{A useful lemma for non-simple $E$.}

Our analysis of the dominating irreducible components of the variety
of pairs $(C,A)$, where $C$ is a smooth curve in $|2H|$, $A$ a
pencil on $C$, and the Lazarsfeld-Mukai associated vector bundle
$E(C,A)$ is not simple, is based on the following key lemma (see
\cite{GL87}, \cite{DM89}. See \cite[Lemma 2.1]{CP95} for this
precise statement).

\begin{lem}\label{CPextension}
Let $S$ be a $K3$ surface, $C$ a smooth curve on $S$ and $A$ a
base-point-free line bundle on $C$, such that $h^0(C,A)=2$. If
$E(C,A)$ is not a simple vector bundle, then there exists two line
bundles $M$ and $N$ on $S$ and a $0$-dimensional subscheme $\xi$ of
$S$ such that
\begin{enumerate}
 \item $h^0(S,M)\geq 2$, $h^0(S,N)\geq 2$;
 \item $N$ is base-point-free;
 \item there is an exact sequence
 $$
 0\to M \to E \to N\otimes I_{\xi} \to 0.
 $$
\end{enumerate}
Moreover if $h^0(M-N)=0$, then $\supp(\xi)=\emptyset$ and the above
sequence is split.
\end{lem}

This result allows to describe the parameter space for non-simple
Lazarsfeld-Mukai vector bundles $E(C,A)$ as an open subset of a
projective bundle over the Hilbert scheme of points on $S$. We shall
use this lemma for $\lg(\xi)=2$ or $\xi=\emptyset$.

\section{Pencils on curves in $|2H|$.}

\subsection{The invariants of curves in $|2H|$.}

Let $C$ be a smooth curve in the linear system. The adjunction
formula computes the genus $g(C)=4r-3$.

The other basic invariants, the Clifford index and the gonality are
obtained as follows.

Observe that $H|_C$ contributes to the Clifford index
of $C$, and
$$
\Cliff(H|_C)=2r-4.
$$
In particular,
$\Cliff(C)\le 2r-4$, whereas the Clifford index
of generic curves of genus $4r-3$ equals $2r-2$. By
the main result of \cite{GL87}, the Clifford index
of $C$ is computed by a line bundle $N=aH+b\ell$ on
the surface $S$. Let $M=(2-a)H-b\ell$.
The bundle $N$ can be assumed to be base-point-free,
and $h^1(N)=0$, and $h^0(N)=h^0(N|_C)$, and
we may assume that $h^1(M)=0$, and the restriction
map $H^0(M)\to H^0(M|_C)$ is surjective (compare to
\cite[Proof of Proposition 3.3]{CP95}). Since both
$N|_C$, and $M|_C$ compute the Clifford index of $C$, we
obtain $a=1$. The other conditions imply $b=0$, hence
$$
\Cliff(C)=2r-4.
$$

Finally, apply \cite[Proposition 3.3]{CP95} to conclude the
the gonality of $C$ equals $2r-2$.

\subsection{The parameter space of pairs, and the proof strategy.}

The main ingredient used in the proof of Theorem \ref{2h} is the
parameter space of pairs $(C,A)$, where $C\in |2H|$ is a smooth
curve, and $A$ is a pencil of given degree $m$ on $C$. Denote
$|2H|_s$ the open subset of $|2H|$ corresponding to smooth curves.
Following \cite{AC81}, there exists a variety
$\mathcal{W}^1_m(|2H|_s)$, and a projective morphism
$$
\pi_S:\mathcal{W}^1_m(|2H|_s)\to |2H|_s,
$$
whose fibre $\pi_S^{-1}[C]$ over any $C\in |2H|_s$ is
scheme-theoretically isomorphic to the variety of special divisors
$W^1_m(C)$. The proof idea of Theorem \ref{2h} is to estimate the
relative dimension of the dominating irreducible components of
$\mathcal{W}^1_{2r-2+n}(|2H|_s)$, for $n\in\{0,1,2\}$.

Notice that the case $n=0$ has already been settled in \cite[Theorem
3.1, and Lemma 3.2]{CP95}: a generic member of the linear system
$|2H|$ carries only finitely many minimal pencils (remark that
$\rho(4r-3, 1, 2r-2)=-3$).

\medskip

The case $n=1$ will follow as an immediate consequence of
\cite[Theorem 3.1, and Lemma 3.2]{CP95}, and of a Lemma of Accola
(see \cite{Ac81}, cf. \cite[Lemma 3.1]{ELMS89}, and \cite{CM91}).

\begin{lem}\label{n=1}
There is no base-point-free $\mathfrak{g}^1_{2r-1}$ on a smooth $C\in |2H|$.
\end{lem}

\begin{proof}
Indeed, if $A$ was such a pencil, then by Accola's Lemma, we would
have (the bundle $\sO_C(H)$ is semi-canonical)
$$
 2h^0(C,\sO_C(H-A))\ge 2 h^0(C,\sO_C(H))-(2r-1).
$$

Since $h^0(S,\sO_S(H))=h^0(C,\sO_C(H))=r+1$, we obtain
$2h^0(C,\sO_C(H-A))\ge 3$, i.e. $h^0(C,\sO_C(H-A))\ge 2$. In this
case, $\sO_C(H-A)$ would be a pencil of degree $2r-3$, contradicting
$\gon(C)=2r-2$.
\end{proof}

Therefore,  for smooth curves $C\in |2H|$, the variety
$W^1_{2r-1}(C)$ is $1$-dimensional, and its irreducible components
are  obtained by adding base-points to the (finitely many)
$\mathfrak{g}^1_{2r-2}$'s on $C$.

\medskip

So, we may suppose $n=2$. Let ${\mathcal{W}}$ be an irreducible
component of $\mathcal{W}^1_{2r}(|2H|_s)$, such that for a generic
pair $(C,A)$, the line bundle $A$ is base-point-free. By \cite[Lemma
3.5, p. 182]{ACGH85}, we necessarily have $h^0(C,A)=2$. We
distinguish two cases according to the behaviour of the associated
Lazarsfeld-Mukai bundle $E(C,A)$. Precisely, components whose
generic member $(C,A)$ has {\it simple}, respectively {\it non
simple}, Lazarsfeld-Mukai vector bundle are called {\it simple
components}, respectively {\it non simple components}.

These two cases are treated in separate Subsections in the sequel.

\subsection{The study of simple components.}

In this Subsection, we study components of $\mathcal{W}^1_{2r-2+n}$
whose generic member  $(C,A)$ gives rise to a {\it simple}
Lazarsfeld-Mukai vector bundle $E(C,A)$.  The relative dimension of
{\it simple components} dominating $|2H|$ is determined thanks to a
more general result due to Pareschi.

\begin{thm}[{\cite[Theorem 2, p. 196]{P95}}]
\label{Pareschi} Let $S$ be a $K3$ surface, and $L$ a line bundle on
$S$. Let $\mathcal{W}$ be an irreducible component of
$\mathcal{W}^r_d(|L|_s)$. Suppose that for a generic pair $(C,A)\in
\mathcal{W}$, the line bundle $A$ is base-point-free, $h^0(C,A)=r+1$
and the associated vector bundle $E(C,A)$ is simple. If the Petri
map $\mu_{0,A}$ is not injective, then the differential
$$
 (d\pi_S)_{|(C,A)} : T_{(C,A)} \mathcal{W} \to T_C |L|=H^0(C,K_C)
$$
is not surjective.
\end{thm}

By standard Brill-Noether theory \cite{ACGH85}
this result immediately yields the following.

\begin{cor}\label{simpledim}
 The relative dimension of the simple components of  $\mathcal{W}^r_d(|L|_s)$
 dominating $|L|$ equals $\rho(p_a(L), r, d)$.
\end{cor}

Notice that in our case, $\rho(4r-3, 1, 2r-2+n)$ is negative for
$n\in\{0,1\}$, it equals $1$, for $n=2$, and it equals $3$, for
$n=3$. In particular, there is no dominating simple component in the
cases $n=0$ and $n=1$.

\medskip

Theorem \ref{Pareschi} is stated, and proved, in \cite{P95} under
the hypothesis that the linear system $|L|$ does not contain
reducible or multiple curves. This condition on $|L|$ is only used
in order to insure the simplicity of the associated bundles
$E(C,A)$, which we assume. For the convenience of the reader, since
we have slightly changed the hypothesis, we sketch below Pareschi's
very nice infinitesimal argument.

\smallskip

\noindent
{\it Proof of Theorem \ref{Pareschi}.}
Consider the map
$$
 \mu_{1,A,S} : \ker(\mu_{0,A})\to H^1(C, \sO_C)
$$
defined as the composition of the Gaussian map (see \cite{AC81})
$$
 \mu_{1,A} :\ker(\mu_{0,A})\to H^0(C, 2K_C),
$$
with the transpose $\delta^{\vee}_{C,S}$ of the Kodaira-Spencer map
$$
 \delta_{C,S} : H^0(C,N_{C/S}=K_C) \to H^1(C, T_C).
$$

The first fact is that, by standard first-order deformation theory,
(see for instance [CGGH, \S2(c)]) one has
\begin{equation}\label{ann}
 \Im (d\pi_S)_{|(C,A)} \subset \Ann (\mu_{1,A,S})\subset
 H^1(C,\sO_C)^*\cong H^0(C,K_C).
\end{equation}

On the other hand, Pareschi \cite[Lemma 1]{P95} shows that, up to a
scalar factor, $\mu_{1,A,S}$ coincides with the coboundary map:

\begin{equation}\label{coboundary}
 H^0 (C, M_A\otimes K_C(-A)) \to H^1(C,\sO_C)
\end{equation}
of the exact sequence (\ref{FA&MA}) twisted by  $K_C (-A)$ :
$$
 0\to \sO_C \to E(C,A)^* \otimes K_C (-A)\to M_A\otimes K_C (-A)\to 0.
$$

Now, by hypothesis, $\ker(\mu_{0,A})\cong H^0 (C, M_A\otimes
K_C(-A))$ is not zero. Moreover, again by hypothesis,
$$
 h^0(S, E(C,A)\otimes E(C,A)^*)=1.
$$

In particular, since twisting (\ref{E}) with $E(C,A)^*$ we have that
$$
 H^0(S, E(C,A)\otimes E(C,A)^*)\cong H^0(C,E(C,A)^*\otimes K_C(-A)),
$$
by (\ref{coboundary}) we get that $\mu_{1,A,S}$ is injective, hence
$$
\mu_{1,A,S}\ \textrm{is not the zero map.}
$$

Therefore
$$
  \Ann (\mu_{1,A,S}) \subsetneq H^1(C,\sO_C)^*\cong H^0(C,K_C)
$$
and, by (\ref{ann}), the differential ${d\pi_S}_{|(C,A)}$ cannot be
surjective.

\hfill$\Box$

\subsection{The study of non simple components.}

In this section we want to study the irreducible components of
${\mathcal{W}}^1_{2r}(|2H|_s)$,  whose generic
member $(C,A)$ has the property that $A$ is a complete
base-point-free pencil of degree $2r$, with the associated
bundle $E(C,A)$ not simple.

We prove the following.

\begin{thm}\label{nonsimplecomp}
Let ${\mathcal{W}}$ be a non simple irreducible component. If
${\mathcal{W}}$ dominates $|2H|$, then its relative dimension equals
$1$.
\end{thm}

A key tool to study non simple $E:=E(C,A)$ is provided by Lemma
\ref{CPextension}. Accordingly, we have an exact sequence :
\begin{equation}\label{eqn:extension}
 0\to M \to E \to N\otimes I_{\xi} \to 0.
\end{equation}

Since in our case $\Pic(S)=\mathbb{Z}[H]\oplus \mathbb{Z}[\ell]$,
the line bundles $M$ and $N$ are of the form
$$
 M=aH+b\ell,\ \ N=a'H+b' \ell.
$$

As $c_1(E)=2H$, we have that $a+a'=2$, and $b'=-b$, and using the fact that
$h^0(M),h^0(N)\geq 2$, we get
$$
a=a'=1.
$$

Moreover,
$$
 c_2(E)=M.N+\lg(\xi)=2r.
$$

In conclusion, two cases may occur :

\begin{enumerate}
\item[($\star$)] $b=0$ and $\lg(\xi)=2$;
\item[($\star\star$)] $b=1,\ b'=-1$ and $\xi$ is
empty;
\end{enumerate}
the case $b=-1,\ b'=1$ is excluded since $N$ is globally generated.

\medskip

We point out the following useful fact.

\begin{lem}\label{lem:locally free2}
Let $\sE$ be a torsion-free sheaf given by a non-trivial extension
\begin{equation}\label{eqn:torsionfree2}
0\to\sO_S(H)\to\sE\to\sO_S(H)\otimes I_\xi\to 0,
\end{equation}
where $\xi\subset S$ is a zero-dimensional subscheme with
$\lg(\xi)=2$. Then $\sE$ is locally free.
\end{lem}

\begin{proof} Suppose $\sE$ was not locally free, and denote
$E=\sE^{**}$. As in \cite{OSS80}, cf. \cite[Proof of Proposition
3.9]{La97}, there exists a subscheme $\eta\subsetneq\xi$, and a
commutative diagram:
\begin{equation}\label{eqn:eta}
\xymatrix{ & &0\ar[d] &0\ar[d]\\
0\ar[r] &\mathcal{O}_S(H)\ar@{=}[d]\ar[r] &\sE \ar[d]\ar[r]
& \mathcal{O}_S(H)\otimes I_\xi \ar[d]\ar[r] &0\\
0\ar[r] &\mathcal{O}_S(H)\ar[r] &E\ar[d]\ar[r] &
\mathcal{O}_S(H)\otimes I_\eta\ar[d]\ar[r] & 0\\
 & &E/\sE\ar@{=}[r]\ar[d] &I_{\eta\subset\xi}\ar[d]\\
& &0&0}
\end{equation}
where $E/\sE$ is supported on the singular locus of $\sE$. Since
$\lg(\eta)\le 1$, and $E$ is locally free, we deduce that
$\eta=\emptyset$. Then $E$ is split, and therefore the extension
(\ref{eqn:torsionfree2}) is trivial contradicting the hypothesis.
\end{proof}

\begin{lem}\label{lem:unique}
For any non simple bundle $E$, the extension (\ref{eqn:extension})
is uniquely determined.
\end{lem}

\begin{proof} We observe first that $h^0(S,E(-M))=1$. Indeed, in the
case ($\star$), we use $h^0(S,I_\xi)=0$, as $\xi\ne\emptyset$. In
the case ($\star\star$), we use $h^0(S,\sO_S(-2\ell))=0$.

Next, we prove that one $E$ cannot lie in two different extensions:
$$
0\to\sO_S(H)\to E \to\sO_S(H)\otimes I_{\xi} \to 0,
$$
and
$$
0\to\sO_S(H+\ell)\to E \to\sO_S(H-\ell)\to 0.
$$

If we found two extensions as above, then, using $h^0(S,E(-H))=1$ as
remarked above, we would be led to a commutative diagram:
\begin{equation}
\xymatrix{ \sO_S(H)\ar@^{(->}[r]\ar@^{(->}[d]&E\ar@{=}[d]\\
\sO_S(H+\ell)\ar@^{(->}[r] & E}
\end{equation}
Furthermore, by obvious reasons, we would obtain an inclusion
$$
\sO_S(H+\ell)/\sO_S(H)\hookrightarrow E/\sO_S(H).
$$
This is absurd, as $E/\sO_S(H)=\sO_S(H)\otimes I_{\xi}$ is
torsion-free, whereas $\sO_S(H+\ell)/\sO_S(H)$ is supported on
$\ell$.
\end{proof}

Then, since the data of an extension (\ref{eqn:extension}) consists
of a (possibly empty) $0$-dimensional subscheme $\xi$ of S and an
element $t\in \mathbb{P}H^1(S,M^{\vee}\otimes N\otimes I_{\xi})^*$,
we obtain:

\begin{prop} The parameter space $\mathcal{P}$ for the non simple vector
bundles $E$ is either isomorphic to $S^{[2]}$ in the case ($\star$),
or to a projective plane $\mathbb{P}^2$ in the case ($\star\star$).
\end{prop}

\begin{proof} In the case $(\star)$, the parameter space in question
parametrises pairs $(t,\xi)$ with $\xi \in S^{[2]}$, and $t\in
\mathbb{P}H^1(S,M^{\vee}\otimes N\otimes I_{\xi})^* =
\mathbb{P}H^1(S,I_{\xi})^*$ {\em arbitrary}. Since
$h^1(S,I_{\xi})=1$, the parameter space  coincides actually with
$S^{[2]}$, {\em see} also Lemma \ref{lem:locally free2}. In the case
$(\star\star)$, we observe that the split bundle is not of type
$E(C,A)$, as the union of the zero loci of two sections: one of
$\mathcal{O}_S(H+\ell)$ and the other one of $\mathcal{O}_S(H-\ell)$
is a reducible curve in $|2H|$.
\end{proof}

\begin{rmk}
{\rm The existence of a universal extension
over $S^{[2]}$ is insured by a more general
result due to H. Lange (see
\cite[Proposition 4.2 and Remark 3.5]{Lange}).}
\end{rmk}

Now consider the Grassmann bundle
$\mathcal{G}\buildrel{p}\over{\to}\mathcal{P}$, whose
fiber at a point $[E]\in \mathcal{P}$ is the Grassmannian of
two-dimensional subspaces of the global sections of $E$:
$$
 p^{-1} ([E]) = \Gr(2, H^0(S,E)).
$$

The dimension of $\mathcal{G}$ is
\begin{equation}\label{dimG}
 \dim (\mathcal{G}) = \dim(\mathcal{P}) +
 \dim( \Gr(2, H^0(S,E))) = \dim \mathcal{P} + 4r-4
\end{equation}
and there is a rational map
\begin{equation}\label{d}
 d : \mathcal{G} \dashrightarrow |2H|,\ \ \ ([E], \Lambda ) \mapsto d_E (\Lambda)
\end{equation}
where the map $d_E$ is the determinant map, defined in (\ref{d_E}).
The utility of $\mathcal{G} \buildrel{d}\over{\dashrightarrow}
|2H|$ is made clear by the following result.

\bigskip

\begin{prop}\label{prop:G&W}
The irreducible components
${\mathcal{W}}$ of ${\mathcal{W}}^1_{2r}(|2H|_s)$ whose generic
member $(C,A)$ is such that $A$ is base-point-free and the
associated bundle $E(C,A)$ is not simple, are birational to the two
Grassmann bundles $\mathcal G_{\star}$ and $\mathcal G_{\star\star}$
corresponding to the cases $(\star)$ and $(\star\star)$. Moreover,
denoting these components by ${\mathcal{W}}_\star$ and
${\mathcal{W}}_{\star\star}$, and by $f_\star$ and $f_{\star\star}$
the birational maps between them and the Grassmann bundles, we have
that $f_\star$ and $f_{\star\star}$ commute with the maps $d$
and $\pi_S$ on $|2H|$.
\end{prop}

\begin{proof}
We write the proof only for the case ($\star$). The map $f_\star$
from $\mathcal{G}_\star$ to the corresponding irreducible component
${\mathcal{W}_\star}$ associates to a generic pair $([E], \Lambda)$
in the Grassmann bundle, the element $(C,A)\in {\mathcal{W}_\star}$,
where $C:= d_E(\Lambda)$ and
 $$
  A:= \Im (\Lambda\otimes \sO_C \hookrightarrow E\otimes \sO_C).
 $$

The map $f_\star$ has degree one, since a pair $(C,A)$, with $C$
smooth, and $A$ base-point-free with $h^0(C,A)=2$, determines a
unique element $([E(C,A)], H^0(C,A)^*)\in \mathcal{G}_\star$. The
commutativity of the diagram
\begin{equation}\label{G&Wdiag}
\xymatrix{\mathcal{W}_\star \ar[dr]
\ar[dr]_{(\pi_S)_{|\mathcal{W}_\star}} &
\mathcal{G}_\star\ar@{-->}[l]^{f_\star}\ar@{-->}[d]^{d_{|\mathcal{G}_\star}}\\
 &
|2H|\\
}
\end{equation}
follows immediately from the description of the birational map
$f_\star$ we have given.
\end{proof}

In case ($\star\star$) the dimension of $\mathcal P$ equals $2$, so we get
$$
 \dim (\mathcal G_{\star\star})= 4r-2.
$$

Since $\dim(|2H|)=4r-3$, if $d$ is dominant, then
\begin{equation}\label{dimG_**}
 \dim (d^{-1}(C))=1,
\end{equation}
implying that the relative dimension of the component
$\mathcal{W}_{\star\star}$ equals one.

\medskip

The case $(\star)$ is slightly different. Notice that in this case,
by (\ref{dimG}), we have
\begin{equation}\label{dimG_*}
 \dim (\mathcal{G}_\star)  = 4r.
\end{equation}
So if $d:\mathcal{G}_\star\dashrightarrow |2H|$ were dominant, the
varieties of pencils $W^1_{2r} (C)$ of a generic curve $C\in |2H|$
would have dimension equal to $3$, and thus they would not satisfy
the linear growth conditions. Theorem \ref{nonsimplecomp} will then
be a consequence of the following result.
\begin{lem}\label{g*notdom}
 The Grassmann bundle $\mathcal{G}_\star$ does not dominate $|2H|.$
\end{lem}

\begin{proof} We use again Pareschi's infinitesimal approach. Suppose
that $\mathcal{W}_\star$ dominates $|2H|$. If $(C,A)\in
\mathcal{W}_\star$ is a generic pair, then arguing as in the proof
of Theorem \ref{Pareschi} one obtains that $\ker (\mu_{0,A})$,
which, by the base-point-free pencil trick is isomorphic to
$H^0(C,K_C(-2A))$, is at least two-dimensional.

On the other hand, if $\mathcal{W}_\star$ dominates $|2H|$,
then (\ref{ann}) implies that $\mu_{1,A,S}\equiv 0$, hence
we have an exact sequence
$$
0 \to H^0(C,\mathcal{O}_C)\to H^0(C,E^*|_C\otimes K_C(-A)) \to
H^0(C,K_C(-2A))\to 0.
$$

We would obtain then $h^0(C,E^*|_C\otimes K_C(-A))\ge 3$.

On the other hand, twisting the exact sequence
$$
0 \to \mathcal{O}_S(H)\to E \to
\mathcal{O}_S(H)\otimes I_\xi \to 0
$$
by $E^*$, and recalling that by determinant reasons,
$E^*\cong E\otimes \mathcal{O}_S(-2H)$, we get
$$
0\to E(-H)\to E\otimes E^*\to E(-H)\otimes I_\xi\to 0,
$$
implying $h^0(S,E\otimes E^*)\le 2$. Furthermore, the exact sequence
$$
0\to \Lambda\otimes \mathcal{O}_S\to E\to K_C(-A)\to 0
$$
twisted by $E^*$ together with the relations $h^0(E^*)=h^1(E^*)=0$,
yields to an isomorphism
$$
H^0(S,E\otimes E^*)\cong H^0(C,E^*|_C\otimes K_C(-A)),
$$
which lead to a contradiction.

Consequently, the map $\mu_{1,A,S}$ is not identically zero and the
irreducible component $\mathcal W$ of
$\mathcal{W}^1_{2r}(|2H|_s)$ cannot be dominant by (\ref{ann}).
\end{proof}

We pass now to the proof of Theorem \ref{nonsimplecomp}.
\begin{proof}[Proof of Theorem \ref{nonsimplecomp}]
Non simple components of ${\mathcal{W}}^1_{2r}(|2H|_s)$,
are birational to one of the two Grassmann bundles
$\mathcal{G}_\star$ and  $\mathcal{G}_{\star\star}$, see Proposition \ref{prop:G&W}.

The first one cannot dominate $|2H|$, by Lemma \ref{g*notdom}.
If the latter one dominates $|2H|$, then its relative dimension equals
$1$, by (\ref{dimG_**}), which we wanted to prove.
\end{proof}

\subsection{Conclusion of proofs of Theorem \ref{2h}
and Proposition \ref{counter-example}.}

As we have analysed several cases, we now make the point and put
them together to show how they imply Theorem \ref{2h}. Also, we
prove that curves in the image of the non-simple component
birational to $\mathcal{G}_\star$ violate Statement (T).

\begin{proof}[Proof of Theorem \ref{2h}]
Let $C$ be a generic curve in $|2H|$. Let $W^1_{2r-2+n}(C),$
$n\in\{0,1,2\}$ be the variety of degree $2r-2+n$ pencils on C.  For
$n=0$, by \cite{CP95} the dimension of  $W^1_{2r-2}(C)$ is zero. For
$n=1$, by Lemma \ref{n=1},  $W^1_{2r-1}(C)$ has only irreducible
components $W$ whose member is a pencil $A'$ on $C$ obtained by
adding a base point to a $A\in W^1_{2r-2}(C)$  Hence $\dim
(W^1_{2r-1}(C))=1$. For $n=2$, and $W$ an irreducible
component of $W^1_{2r}(C)$, we have several possibilities.
\begin{enumerate}
\item[\bf{(a)}] A generic $A\in W$ is base-point-free.
\item[\bf{(b)}] Any $A\in W$ has base-points.
\end{enumerate}

In case {\bf(a)}, also using \cite[Lemma 3.5, p. 182]{ACGH85}, we
obtain two subcases.
\begin{enumerate}
\item[\bf{(a1)}] A generic $A\in W$ is such that $h^0(C,A)=2$ and $E(C,A)$ is simple.
\item[\bf{(a2)}] A generic $A\in W$ is such that $h^0(C,A)=2$ and $E(C,A)$ is not simple.
\end{enumerate}

In case {\bf(a1)}, the dimension of such a component is given by
Corollary \ref{simpledim}, and equals $1$, if $n=2$.

In case {\bf(a2)}, as $C$ is generic, and, by Lemma \ref{g*notdom},
$\mathcal{G}_\star$ does not dominate $|2H|$, we have that $W$ must be
birational to the fibre of $\mathcal{G}_{\star\star}$ over $C$, by
Proposition \ref{prop:G&W}. So its dimension equals $1$, by
(\ref{dimG_**}).

In case {\bf(b)}, the component $W$ is dominated by an irreducible
component of $W^1_{2r-1}(C)\times C$, so we are done, using case
{\bf(a)}, and $n=0,1$.

The conclusion is that for $n=0,1,2$ we have
$$
 \dim (W^1_{2r-2+n}(C)) =n,
$$
for $C$ generic in $|2H|$, and the theorem is proved.
\end{proof}

\begin{proof}[Proof of Proposition \ref{counter-example}] For $n=0$,
as thanks to \cite{CP95} the gonality of any smooth curve in $|2H|$ is
$2r-2$, by \cite{FHL84} we necessarily have
$$
  \dim(W^1_{2r-2}(C)) \leq 1.
$$

For $n=1$, thanks to Accola's Lemma \cite{Ac81}, for any smooth
curve $C$ in $|2H|$ we have
$$
  \dim(W^1_{2r-1}(C))=\dim(W^1_{2r-2}(C)) +1\le 2.
$$

For $n=2$, consider the irreducible component of
$\mathcal{W}^1_{2r}(|2H|_s)$ which is birational to
$\mathcal{G}_\star$ (notice that $\mathcal{G}_\star$ is not empty by
the Serre construction). As we have proved in the previous
subsection,  $\mathcal{G}_\star$ does not dominate $|2H|$. Hence,
taking {\it any} smooth curve $C$ lying in the image of
$\mathcal{G}_\star$, and using (\ref{dimG_*}) and \cite{FHL84} as we
have done before, we get $\dim(W^1_{2r}(C))\ge 4$, and
$\dim(W^1_{2r+1}(C))\ge 5$. By \cite{FHL84}, we obtain that the
dimension of $W^1_{2r-1}(C)$ is at least two; in particular, it
equals two. It implies $\dim(W^1_{2r-2}(C))=1$, and
$\dim(W^1_{2r}(C))=4$.
\end{proof}

\section{Green's conjecture for curves in $|2H+\ell|$.}

\subsection{Passing from $|2H|$ to $|2H+\ell|$.}

We recall first Green's Hyperplane Section Theorem \cite[Theorem
(3.b.7)]{G84}, which reads, in our case
$$
K_{p,1}(S,\sO_S(2H))\cong K_{p,1}(X,\sO_X(2H)),
$$
for any connected curve $X\in |2H+\ell|$, and any positive integer
$p$. We apply this result twice, once for a smooth curve, and once
again for a curve with two reducible components. Specifically, we
obtain $K_{p,1}(X,K_X)\cong K_{p,1}(C+\ell,\omega_{C+\ell})$, where
$X\in |2H+\ell|$, and $C\in |2H|$ are smooth curves, and
$C.\ell=x+y$. Next, remark that the exact sequence:
$$
0\to \sO_\ell(-2)\to \sO_{C+\ell}\to \sO_C\to 0
$$
yields, after tensoring with $\sO_S(2H+\ell)$, to an isomorphism of
vector spaces $H^0(C+\ell,\omega_{C+\ell})\cong H^0(C,K_C(x+y))$.
Similarly, we obtain an inclusion
$H^0(C+\ell,\omega_{C+\ell}^{\otimes 2})\subset
H^0(C,(K_C(x+y))^{\otimes 2})$. By the definition of the Koszul
cohomology groups, we get an isomorphism
$$
K_{p,1}(C+\ell,\omega_{C+\ell})\cong K_{p,1}(C,K_C(x+y)).
$$

The genus of $X$ is $4r-2$, its gonality equals $2r$, and its
Clifford index equals $2r-3$, \cite{ELMS89}. Green's conjecture for
$X$ predicts
$$
K_{2r,1}(X,K_X)=0.
$$

It amounts to prove (by what we have said above):
$$
K_{2r,1}(C,K_C(x+y))=0.
$$

The curve $C$ is of genus $4r-3$, gonality $2r-2$, and Clifford
index $2r-4$, (\cite{ELMS89}, \cite{CP95}, and \cite{GL87}). Note
that $h^0(C,K_C(x+y))=4r-2$, and the vanishing
$K_{2r,1}(C,K_C(x+y))=0$ is the one predicted by the
Green-Lazarsfeld conjecture for the bundle $K_C(x+y)$.

\subsection{Pencils through $x+y$.}

We prove the following.
\begin{lem}\label{ThroughXY}
Let $C\in |2H|$ be any smooth curve, and $x,y\in  C$ its
intersection points with the line $\ell$. For any integer $n\ge 0$,
there is no base-point-free line bundle $A$ on $C$ with
$h^0(C,A)=2$, $\deg(A)=2r-2+n$, and $h^0(C,A(-x-y))\ne 0$.
\end{lem}
\begin{proof}
We argue by contradiction. Suppose there exists complete
base-point-free pencil $A$ of degree $(2r-2+n)$ on $C$ such that
$$
 h^0(C,A(-x-y))\not= 0.
$$

Since $A$ is base-point-free, and $h^0(C,A)=2$, we necessarily have
$$
 h^0(C,A(-x-y)) =1.
$$

Consider the associated vector bundle $E:=E(C,A)$.
Twisting the exact sequence (\ref{E}) by $\mathcal{O}_S(\ell)$
we obtain
\begin{equation}\label{El}
 0\to H^0(C,A)^*\otimes \sO_S(\ell)
 \to E\otimes \mathcal{O}_S(\ell)\to K_C(-A+x+y)\to 0.
\end{equation}
Using the Riemann-Roch theorem for surfaces, one checks that
$$
h^1 (S, \sO_S(\ell))=0.
$$
Therefore from (\ref{El}) we get
$$
 h^0 (S,E\otimes \mathcal{O}_S(\ell))=
 2+ h^0(K_C(-A+x+y))= 2+h^1(A(-x-y))= 2r+3-n,
$$
where the last equality also follows from Riemann-Roch.

On the other hand from (\ref{E}) one also gets
$$
 h^0(S,E)= 2+h^1(C,A)= 2r+2-n.
$$
Consider the exact sequence defining $\ell$, twisted by
$E\otimes\mathcal{O}_S (\ell)$ :
$$
 0\to E\to E\otimes\mathcal{O}_S (\ell) \to E_{|\ell} (-2) \to 0.
$$
As $h^1(S,E)=0$, from the above exact sequence we get
\begin{equation}\label{=1}
 h^0(E_{|\ell} (-2))= h^0(E\otimes\mathcal{O}_S (\ell)) -
 h^0(E)
 =1.
\end{equation}
Now
$$
 E_{|\ell} (-2)= \mathcal{O}_{\ell} (a) \oplus \mathcal{O}_{\ell} (b).
$$
Since $c_1(E)=2H$, we have $a+b=0$.
So we can rewrite (\ref{=1}) :
$$
 h^0(\ell, \mathcal{O}_{\ell} (a))+h^0(\ell, \mathcal{O}_{\ell} (-a))=1
$$
which is absurd.
\end{proof}

\begin{rmk}\label{rmk:n=0}
{\rm In the case $n=0$, as the gonality of any smooth curve in
$|2H|$ equals $2r-2$, from the previous Lemma we deduce that for any
smooth curve $C$ in the linear system $|2H|$, there is no
$\mathfrak{g}^1_{2r-2}$ passing through the intersection points of
$C$ with the line $\ell$.}
\end{rmk}

We arrive which makes crucial use of Theorem \ref{2h}.

\begin{prop}\label{nolb}
Let $C$ be a generic curve in the linear system $|2H|$, and
$\{x_0,y_0\}=C\cap\ell$. For three generic cycles
$x_1+y_1,x_2+y_2,x_3+y_3\in C^{(2)}$, and for any $n\in\{1,2,3\}$,
there is no line bundle $A\in W^1_{2r-2+n}(C)$, verifying
$$
h^0(C,A(-x_0-y_0))\ne 0,
$$
and
$$
h^0(C,A(-x_{i_j}-y_{i_j}))\ne 0
$$
for any set of indices
$\{i_1,\dots i_n\}\subset \{1,2,3\}$.
\end{prop}

\begin{proof} We argue as in \cite{A05} with the difference that we cannot
simply invoke the genericity of the pairs of points, since $x_0+y_0$
is fixed. Nevertheless, we can overpass this particularity using
Lemma \ref{ThroughXY}.

We have the following.
\medskip
\begin{claim} The incidence variety inside
$\mathop\prod\limits_{n}C^{(2)}\times W^1_{2r-2+n}(C)$,
$$
\Xi:=\{(x_1+y_1,\dots,x_n+y_n,A),
\, h^0(C,A(-x_i-y_i))\ge 1, \mbox{ for all } i=0,\dots,n\}
$$
is at most $(2n-1)$-dimensional.
\end{claim}
\medskip
The proof of the claim is done by analysing
all the possible irreducible components $\mathcal{W}$ of the
universal family $\mathcal{W}^1_{2r-2+n}(|2H|_s)$ dominating $|2H|$,
and applying Lemma \ref{ThroughXY}. We have two cases according
to the behaviour of a generic point $(C,A)\in\mathcal{W}$.
\begin{enumerate}
\item[\bf{(a)}] $A$ is base-point-free.
\item[\bf{(b)}] $A$ has base-points.
\end{enumerate}
In case {\bf(a)}, by \cite[Lemma 3.5, p.182]{ACGH85} we necessarily
have $h^0(C,A)=2$, and hence we may apply Lemma \ref{ThroughXY} to
get to a contradiction. In conclusion this case does not occur.

In case {\bf(b)}, let $z\in C$ be a base-point of $A$. Up to
subtracting from $A$ all the base-points different from $x_0$ and
$y_0$ and applying Lemma  \ref{ThroughXY} to get to a contradiction,
we may assume $z\in \{x_0, y_0\}$.  We consider the incidence
variety $\Xi'$ inside $\mathop\prod\limits_{n}C^{(2)}\times
W^1_{2r-3+n}(C)$ :
 $$
\Xi':=\{(x_1+y_1,\dots,x_n+y_n,A'),
\, h^0(C,A'(-x_i-y_i))\ge 1, \mbox{ for all } i=1,\dots,n\}
$$
(notice that the difference between $\Xi$ and $\Xi'$ is that for the
latter we are not imposing that the line bundles $A'$ pass through
$x_0+y_0$). We have two injective maps from $\Xi'$ to $\Xi$ :
$$
 j_{x_o} : \Xi' \hookrightarrow \Xi ; \ \ A'\mapsto A:= A'+x_0
$$
and
$$
 j_{y_o} : \Xi' \hookrightarrow \Xi ; \ \ A'\mapsto A:= A'+y_0.
$$
Consider the images $ j_{x_o} (\Xi')$ and $ j_{y_o} (\Xi')$. They
are closed inside $\Xi$, and moreover, thanks to Lemma
\ref{ThroughXY}, we have
$$
 \Xi = j_{x_o} (\Xi')\cup j_{y_o} (\Xi').
$$
In particular, the dimension of $\Xi$ equals that of $\Xi'$.
Then we may argue as in
\cite[p.394]{A05}
 and apply Theorem \ref{2h} to conclude that $\Xi$ is
at most $(2n-1)$-dimensional.

So the claim, and hence the proposition is proved.
\end{proof}

\begin{cor}\label{GLfor2H}
Let $C\in|2H|$ be a generic curve, and $x,y$ the
intersection points with the line $\ell$. Then $K_{2r,1}(C,K_C(x+y))=0$;
in particular, the Green-Lazarsfeld conjecture holds
for $C$.
\end{cor}

\begin{proof}
We choose three generic cycles $x_1+y_1,x_2+y_2,x_3+y_3\in C^{(2)}$,
and denote $Y$ be the nodal curve obtained gluing together $x$ with
$y$, and $x_i$ with $y_i$ for all $i$, and we prove that
$K_{2r,1}(Y,\omega_Y)=0$, arguing similarly to the proof of
\cite[Theorem 2]{A05}. We reproduce the arguments here for the
reader's convenience. Assume by contradiction that
$K_{2r,1}(Y,\omega_Y)\not =0$. Then, by the degenerate version of
the Hirschowitz-Ramanan-Voisin result \cite[Proposition 8]{A05},
there exists a rank one torsion-free sheaf $F$ on $Y$, with $\chi(F)
= 2r-2-g(C) $, and $h^0(F)\geq2$. This sheaf is either a line bundle
or the direct image of a line bundle on a partial desingularization
of $Y$ (which cannot be $C$ itself, otherwise we would contradict
the fact that $\gon(C)=2r-2$). Hence $F=\phi_* L$, where $\phi:Z\to
Y$ is a partial normalization of $(4-n)$ of the $4$ nodes of $Y$.
Let $\psi : C\to Z$ be the normalization of the remaining $n$-nodes.
Then, $\chi(Z,L)=\chi(Y,F)=2r-2-g$, and $\chi(C,\psi^*L)= 2r-1-g+n$,
and the latter implies that $\deg(\psi^*L)=2r-2+n$. As $L$ is a
pencil, for each of the $n$ nodes, there exists a non-zero section
vanishing at it. Hence the pencil $\psi^*L$ would contradict Remark
\ref{rmk:n=0} and Proposition \ref{nolb}. So we have proved that
$K_{2r,1}(Y,\omega_Y) =0$.

On the other hand, thanks to \cite [Lemma 2.3]{AV03}, we have :
$$
 K_{2r,1}(C, K_C(x+y))\subset K_{2r,1}(Y,\omega_Y)
$$
and the predicted vanishing is proved. The last assertion
follows using \cite[Theorem 3]{A02}.
\end{proof}

\begin{proof}[Proof of Theorem \ref{2h+l}]
We have to check that
$$
 K_{2r,1}(X,K_X)=0
$$
for a smooth curve $X$ in $|2H+\ell|.$
As we have shown in \S 4.1, we have the equality
$$
 K_{2r,1}(X,K_X) = K_{2r,1}(C,K_C(x+y))
$$
where $C\in |2H|$. The latter Koszul cohomology group vanishes by
the previous Corollary and the Theorem is proved.
\end{proof}

\smallskip

\small

{\sc Institute of Mathematics "Simion Stoilow" of the Romanian
Academy P.O. Box 1-764, RO-014700 Bucharest -- Romania} {\em E-mail
address:} {\tt Marian.Aprodu@imar.ro}

\medskip

\small

{\sc Institut de Recherches Math\'ematiques Avanc\'ees Universit\'e
Louis Pasteur et CNRS 7 rue Ren\'e Descartes, F-67084 Strasbourg
Cedex -- France} {\em E-mail address:} {\tt
pacienza@math.u-strasbg.fr}


\begin{thebibliography}{ACGH85}

\bibitem[Ac81]{Ac81}
R. D. M. Accola, {\em Plane models for Riemann surfaces admitting
certain half-canonical linear series}, In: Riemann Surfaces and
related topics, Proc. of the 1978 Stony Brook Conference, eds. I.
Kra and B. Maskt, Annals of Math. Studies {\bf 97} Princeton Univ.
Press (1981), 7--20.
\bibitem[A02]{A02}
M. Aprodu, {\em On the vanishing of the higher syzygies of curves},
Math. Z. {\bf 241} (2002), 1--15.
\bibitem[A04]{A04}
\rule{1.5cm}{0.25mm} , {\em Green-Lazarsfeld gonality conjecture for
a generic curve of odd genus}, Int. Math. Res. Not. {\bf 63} (2004),
3409--3416.
\bibitem[A05]{A05}
\rule{1.5cm}{0.25mm} , {\em Remarks on syzygies of $d$-gonal
curves}, Math. Res. Lett. {\bf 12} (2005), 387--400.
\bibitem[AV03]{AV03}
M. Aprodu and C. Voisin, {\em Green-Lazarsfeld's conjecture for
generic curves of large gonality}, C. R. Math. Acad. Sci. Paris {\bf
336} (2003), 335--339.
\bibitem[AC81]{AC81}
E. Arbarello and M. Cornalba, {\em Su una congettura di Petri},
Comm. Math. Helv. {\bf 56} (1981), 1--38.
\bibitem[ACGH85]{ACGH85}
E. Arbarello, M. Cornalba, P. A. Griffiths and J. Harris, {\em
Geometry of Algebraic Curves}, Grundlehren. math. Wiss. {\bf 267}
(1985) Springer Verlag.
\bibitem[CGGH83]{CGGH83} J. Carlson, M. Green, Ph. Griffiths, J. Harris,
{\em Infinitesimal variations of Hodge structure. I}, Compositio
Math.  {\bf 50} (1983), no. 2-3, 109--205.
\bibitem[CP95]{CP95}
C. Ciliberto and G. Pareschi, {\em Pencils of minimal degree on
curves on a $K3$ surface}. J. Reine Angew. Math. {\bf 460} (1995),
15--36.
\bibitem[CKM92]{CKM92}
M. Coppens, Ch. Keem and G. Martens, {\em Primitive linear series on
curves}, Manuscripta Math. {\bf 77} (1992), 237--264.
\bibitem[CM91]{CM91}
M. Coppens and G. Martens, {\em Secant spaces and Clifford's
theorem}, Compositio Math.  {\bf 78} (1991), 193--212.
\bibitem[DM89]{DM89}
R. Donagi and D. Morrison, {\em Linear systems on $K3$ sections}, J.
Diff. Geom. {\bf 29} (1989), 49--64.
\bibitem[ELMS89]{ELMS89}
D. Eisenbud, H. Lange, Herbert, G. Martens and F.-O. Schreyer, {\em
The Clifford dimension of a projective curve}, Compositio Math. {\bf
72} (1989), 173--204.
\bibitem[FHL84]{FHL84}
W. Fulton, J. Harris, R. Lazarsfeld, {\em Excess linear series on an
algebraic curve}, Proc. Amer. Math. Soc.  {\bf 92}  (1984),
320--322.
\bibitem[G84]{G84}
M. Green, {\em Koszul cohomology and the geometry of projective
varieties}, J. Diff. Geom. {\bf 19} (1984), 125--171, with an
Appendix by M. Green and R. Lazarsfeld.
\bibitem[GL85]{GL85}
M. Green and R. Lazarsfeld, {\em On the projective normality of
complete linear series on an algebraic curve}, Invent. Math. {\bf
83} (1985), 73--90.
\bibitem[GL87]{GL87}
\rule{1.5cm}{0.25mm} , {\em Special divisors on curves on a $K3$
surface}, Invent. Math. {\bf 89} (1987), 357--370.
\bibitem[GH78]{GH78}
P. Griffiths and J. Harris, {\em Residues and zero-cycles on
algebraic varieties}, Ann. of Math. {\bf 108} (1978), 461--505.
\bibitem[HR98]{HR98}
A. Hirschowitz and S. Ramanan, {\em New evidence for Green's
conjecture on syzygies of canonical curves}, Ann. Sci. \'Ecole Norm.
Sup. {\bf 31} (1998), 145--152.
\bibitem[Ho82]{Ho82}
R. Horiuchi, {\em Gap orders of meromorphic functions on Riemann
surfaces} J. Reine Angew. Math. {\bf 336} (1982), 213--220.
\bibitem[L83]{Lange}
H. Lange, {\em Universal families of extensions}, J. Algebra {\bf
83} (1983), 101--112.
\bibitem[La86]{La86}
R. Lazarsfeld, {\em Brill-Noether-Petri without degenerations}, J.
Diff. Geom. {\bf 23} (1986), 299--307 .
\bibitem[La89]{La89}
\rule{1.5cm}{0.25mm} , {\em A sampling of vector bundle techniques
in the study of linear series}, M. Cornalba (ed.) et al.,
Proceedings of the first college on Riemann surfaces held in
Trieste, Italy, November 9-December 18, 1987. Teaneck, NJ: World
Scientific Publishing Co. (1989) 500--559.
\bibitem[La97]{La97}
\rule{1.5cm}{0.25mm} , {\em Lectures on Linear Series}, IAS/Park
City Math. Series vol. 8 (1997) 163--219.
\bibitem[Lo89]{Lo89}
F. Loose, {\em On the graded Betti numbers of plane algebraic
curves}, Manuscripta Math.  {\bf 64} (1989), 503--514.
\bibitem[Ma82]{Ma82}
G. Martens, {\em \"Uber den Clifford-Index algebraischer Kurven}, J.
Reine Angew. Math. {\bf 336} (1982) 83--90.
\bibitem[Ma84]{Ma84}
G. Martens, {\em On dimension theorems of the varieties of special
divisors on a curve}, Math. Ann. {\bf 267} (1984), 279--288.
\bibitem[Mu89]{Mu89}
S. Mukai, {\em Biregular classification of Fano $3$-folds and Fano
manifolds of coindex $3$}, Proc. Nat. Ac. Science USA {\bf86}
(1989), 3000--3002.
\bibitem[OSS80]{OSS80}
C. Okonek, M. Schneider, H. Spindler, {\em Vector bundles on complex
projective spaces}, Progress in Math. vol. 3, Birkhaeuser Boston,
1980.
\bibitem[P95]{P95}
G. Pareschi, {\em A proof of Lazarsfeld's Theorem on curves on $K3$
surfaces}, J. Alg. Geom. {\bf 4} (1995), 195--200.
\bibitem[T02]{T02}
M. Teixidor i Bigas, {\em Green's conjecture for the generic
$r$-gonal curve of genus $g\geq 3r-7$}, Duke Math. J. {\bf
111}(2002), 195--222.
\bibitem[V02]{V02}
C. Voisin, {\em Green's generic syzygy conjecture for curves of even
genus lying on a $K3$ surface}, J. Eur. Math. Soc. {\bf 4} (2002),
363--404.
\bibitem[V05]{V05}
\rule{1.5cm}{0.25mm} , {\em Green's canonical syzygy conjecture for
generic curves of odd genus}, Compositio Math. {\bf 141} (2005), no.
5, 1163--1190.
\end{thebibliography}
\end{document}